\documentclass[11pt]{article}
\usepackage{amssymb, amsfonts}
\title{Solvability of linear differential systems in the Liouvillian sense}
\author{R.\,R.\,Gontsov, I.\,V.\,Vyugin\footnote{The author is partially supported by Dynasty Foundation, Simons-IUM fellowship and grant RFBR 12-01-33058.}}
\date{}

\begin{document}
\maketitle
\begin{abstract}
The paper concerns the solvability by quadratures of linear differential systems, which is one 
of the questions of differential Galois theory. We consider systems with regular singular points 
as well as those with (non-resonant) irregular ones and propose some criteria of solvability for
systems whose (formal) exponents are sufficiently small.
\end{abstract}

\section{Introduction}

Consider on the Riemann sphere $\overline{\mathbb C}$ a linear differential system
\begin{eqnarray}\label{syst}
\frac{dy}{dz}=B(z)\,y, \qquad y(z)\in{\mathbb C}^p,
\end{eqnarray}
of $p$ equations with a meromorphic coefficient matrix $B(z)$ having singularities
at points $a_1,\ldots,a_n$. A singular point $z=a_i$ is said to be {\it regular}, if
any solution of the system has at most polynomial growth in any sector of small radius 
with vertex at this point and an opening less than $2\pi$. Otherwise the point 
$z=a_i$ is said to be {\it irregular}.

The {\it Picard--Vessiot extension} of the field ${\mathbb C}(z)$ of rational functions
corresponding to the system (\ref{syst}) is a differential field $F$ obtained by adjoining to
${\mathbb C}(z)$ all entries of a fundamental matrix $Y(z)$ of the system (\ref{syst}). One 
says that the system (\ref{syst}) is {\it solvable by quadratures}, if the entries of the matrix 
$Y(z)$ are expressed in elementary or algebraic functions and their integrals or, more formally,  
if the field $F$ is contained in some extension of ${\mathbb C}(z)$ obtained by adjoining  
algebraic functions, exponentials or integrals: 
$$
{\mathbb C}(z)=F_1\subseteq\ldots\subseteq F_m, \qquad F\subseteq F_m,
$$
where $F_{i+1}=F_i\langle x_i\rangle$ ($i=1,\ldots,m-1$), and either $x_i$ is algebraic over $F_i$, 
or $x_i$ is an exponential of an element in $F_i$, or $x_i$ is an integral of an element in $F_i$. 
Such an extension ${\mathbb C}(z)\subseteq F_m$ is called {\it Liouvillian}, thus solvability by 
quadratures means that the Picard--Vessiot extension $F$ is contained in some Liouvillian extension 
of the field of rational functions.
 
Solvability or non-solvability of a linear differential system by quadratures is related to properties 
of its Galois group. The {\it differential Galois group} $G={\rm Gal}\,(F/{\mathbb C}(z))$ of the system 
(\ref{syst}) (of the Picard--Vessiot extension ${\mathbb C}(z)\subseteq F$) is the group of differential 
automorphisms of the field $F$ (i.~e., automorphisms commuting with differentiation) that preserve elements 
of the field ${\mathbb C}(z)$:
$$ 
G=\Bigl\{\sigma: F\rightarrow F\;\Bigl|\;\sigma\circ\frac d{dz}=\frac d{dz}\circ\sigma,\;\sigma(f)=f\quad
\forall f\in{\mathbb C}(z)\Bigr\}.
$$
As follows from the definition, the image $\sigma(Y)$ of the fundamental matrix $Y(z)$ of the system (\ref{syst})
under any element $\sigma$ of the Galois group is a fundamental matrix of this system again, that is, 
$\sigma(Y)=Y(z)C,\;C\in{\rm GL}(p,{\mathbb C})$. As every element of the Galois group is determined uniquely by 
its action on a fundamental matrix of the system, the Galois group $G$ can be regarded as a subgroup of the matrix
group ${\rm GL}(p,{\mathbb C})$. Moreover this subgroup $G\subset{\rm GL}(p,{\mathbb C})$ is algebraic, i.~e.,
closed in the {\it Zariski topology} of the space ${\rm GL}(p,{\mathbb C})$ (the topology whose closed sets are those
determined by systems of polynomial equations), see \cite[Th. 5.5]{Ka}. 

The Galois group $G$ can be represented as a union of finite number of disjoint connected sets that are open and closed
simultaneously (in the Zariski topology), and the set containing the identity matrix is called the {\it identity component}. 
The identity component $G^0\subset G$ is a normal subgroup of finite index \cite[Lemma 4.5]{Ka}. Due to the Picard--Vessiot 
theorem, solvability of the system (\ref{syst}) by quadratures is equivalent to solvability of the subgroup $G^0$ 
\cite[Th. 5.12]{Ka}, \cite[Ch. 3, Th. 5.1]{Kh}. (Recall that a group $H$ is said to be {\it solvable}, if there exist 
intermediate normal subgroups $\{e\}=H_0\subset H_1\subset\ldots\subset H_m=H$ such that each factor group $H_i/H_{i-1}$ is 
Abelian, $i=1,\ldots,m$.)

Alongside the Galois group one considers the {\it monodromy group} $M$ of the system 
(\ref{syst}) generated by the monodromy matrices $M_1,\ldots,M_n$ corresponding to 
analytic continuation of a fundamental matrix $Y(z)$ around the singular points $a_1,\ldots,a_n$. 
(The matrix $Y(z)$ considered in a neighbourhood of a non-singular point $z_0$ goes to 
$Y(z)M_i$ under an analytic continuation along a simple loop $\gamma_i$ encircled a point 
$a_i$.) As the operation of analytic continuation commutes with differentiation and preserves 
elements of the field ${\mathbb C}(z)$ (single-valued functions), one has $M\subseteq G$. 
Furthermore the Galois group of a system whose singular points are all regular is a closure of 
its monodromy group (in the Zariski topology, see \cite[Ch. 6, Cor. 1.3]{Kh}), hence such a system 
is solvable by quadratures, if and only if the identity component of its monodromy group is solvable. 

We are interested in the cases when the answer to the question concerning solvability of a linear
differential system by quadratures can be given in terms of the coefficient matrix of the 
system. For example, in the case of a Fuchsian system (a particular case of a system with regular
singular points) 
\begin{eqnarray}\label{fuchs}
\frac{dy}{dz}=\biggl(\sum_{i=1}^n\frac{B_i}{z-a_i}\biggr)y, \qquad 
B_i\in{\rm Mat}(p, {\mathbb C}),
\end{eqnarray}
whose coefficients $B_i$ are sufficiently small, Yu.\,S.\,Ilyashenko and A.\,G.\,Khovansky 
\cite{IKh1} have obtained an explicit criterion of solvability. Namely, the following statement 
holds:

{\it There exists $\varepsilon=\varepsilon(n,p)>0$ such that a condition of solvability by
quadratures for the Fuchsian system $(\ref{fuchs})$ with $\|B_i\|<\varepsilon$ takes an explicit
form: the system is solvable by quadratures, if and only if all the matrices $B_i$ are triangular
$($in some basis$)$.}

In this article we propose a refinement of the above assertion in which it is sufficient that the 
eigenvalues of the residue matrices $B_i$ be small (the estimate is given), and we also propose
a generalization to the case of a system with irregular singular points.  

\section{A local form of solutions of a system near its singular points}

A singular point $a_i$ of the system (\ref{syst}) is said to be {\it Fuchsian}, if the coefficient 
matrix $B(z)$ has a simple pole at this point.

Due to Sauvage's theorem, a Fuchsian singular point of a linear differential system is regular
(see \cite[Th. 11.1]{Ha}). However, the coefficient matrix of a system at a regular singular point
may in general have a pole of order greater than one. Let us write the Laurent expansion of the
coefficient matrix $B(z)$ of the system (\ref{syst}) near its singular point $z=a$ in the form
\begin{eqnarray}\label{coeff}
B(z)=\frac{B_{-r-1}}{(z-a)^{r+1}}+\ldots+\frac{B_{-1}}{z-a}+B_0+\ldots,
\qquad B_{-r-1}\ne 0.
\end{eqnarray}
The number $r$  is called the {\it Poincar\'e rank} of the system (\ref{syst}) at this point 
(or the Poincar\'e rank of the singular point $z=a$). For example, the Poincar\'e rank of a 
Fuchsian singularity is equal to zero.
\medskip

The system (\ref{syst}) is said to be {\it Fuchsian}, if all its singular points are Fuchsian 
(then it can be written in the form (\ref{fuchs})). The system whose singular points are all 
regular will be called  {\it regular singular}. 

According to Levelt's theorem \cite{Le}, in a neighbourhood of
each regular singular point $a_i$ of the system (\ref{syst}) there
exists a fundamental matrix of the form
\begin{eqnarray}\label{levelt}
Y_i(z)=U_i(z)(z-a_i)^{A_i}(z-a_i)^{\widetilde E_i},
\end{eqnarray}
where $U_i(z)$ is a holomorphic matrix at the point $a_i$, 
$A_i={\rm diag}(\varphi_i^1,\ldots,\varphi_i^p)$ is
a diagonal integer matrix whose entries $\varphi_i^j$ organize in a
non-increasing sequence, $\widetilde E_i=(1/2\pi{\bf i})\ln\widetilde M_i$ 
is an upper-triangular matrix (the normalized logarithm of the
corresponding monodromy matrix) whose eigenvalues $\rho_i^j$
satisfy the condition
$$
0\leqslant{\rm Re}\,\rho_i^j<1.
$$
Such a fundamental matrix is called a {\it Levelt matrix}, and
one also says that its columns form a {\it Levelt basis} in the
solution space of the system (in a neighborhood of the regular 
singular point $a_i$). The complex numbers $\beta_i^j=\varphi_i^j+\rho_i^j$ 
are called the (Levelt) {\it exponents} of the system at the
regular singular point $a_i$. 

If the singular point $a_i$ is Fuchsian, then the corresponding matrix $U_i(z)$
in the decomposition (\ref{levelt}) is holomorphically invertible at this point,
that is, $\det U_i(a_i)\ne0$. It is not difficult to check that in this case 
the exponents of the system at the point $a_i$ coincide with the eigenvalues of 
the residue matrix $B_i$. And in the general case of a regular singularity $a_i$ 
there are estimates for the order of the function $\det U_i(z)$ at this point 
obtained by E.\,Corel \cite{Co} (see also \cite{Go}):
$$
r_i\leqslant{\rm ord}_{a_i}\det U_i(z)\leqslant
\frac{p(p-1)}2\,r_i,
$$
where $r_i$ is the Poincar\'e rank of the regular singular point $a_i$. These 
estimates imply the inequalities for the sum of exponents of the regular system
over all its singular points, which are called the {\it Fuchs inequalities}:
\begin{eqnarray}\label{fuchsineq}
-\frac{p(p-1)}2\sum_{i=1}^nr_i\leqslant\sum_{i=1}^n\sum_{j=1}^p\beta_i^j
\leqslant-\sum_{i=1}^nr_i
\end{eqnarray}
(the sum of exponents is an integer).

Let us now describe the structure of solutions of the system (\ref{syst}) near one of its irregular singular
points. We suppose that the irregular singularity $a_i$ of Poincar\'e rank $r_i$ is 
{\it non-resonant}, that is, the eigenvalues $b_i^1,\ldots,b_i^p$ of the leading term $B_{-r_i-1}$
of the matrix $B(z)$ in the expansion (\ref{coeff}) at this point are pairwise distinct. Let us
fix a matrix $T_i$ reducing the leading term $B_{-r_i-1}$ to the diagonal form 
$$
T_i^{-1}B_{-r_i-1}T_i={\rm diag}(b_i^1,\ldots,b_i^p).
$$ 
Then the system possesses a uniquely determined formal fundamental matrix $\widehat Y_i(z)$ 
of the form
$$
\widehat Y_i(z)=\widehat F_i(z)(z-a_i)^{\Lambda_i}e^{Q_i(z)},
$$
where
\begin{itemize} 
\item[a)] $\widehat F_i(z)$ is a matrix formal Taylor series in $z-a_i$ and $\widehat F_i(a_i)=T_i$; 

\item[b)] $\Lambda_i$ is a constant diagonal matrix whose diagonal entries are called the
{\it formal exponents} of the system (\ref{syst}) at the irregular singular point $a_i$;

\item[c)] $Q_i(z)={\rm diag}(q_i^1(z),\ldots,q_i^p(z))$ is a diagonal matrix whose diagonal entries 
$q_i^j(z)$ are polynomials in $(z-a_i)^{-1}$ of degree $r_i$ without a constant term, 
$$
q_i^j(z)=-\frac{b_i^j}{r_i}\,(z-a_i)^{-r_i}+o((z-a_i)^{-r_i}). 
$$ 
\end{itemize}

Furthermore a punctured neighbourhood of the point $a_i$ can be covered by a set $\{S_i^1,\ldots,S_i^{N_i}\}$ 
of ''good'' open sectors with vertex at this point (which we take to be arranged in counterclockwise order 
starting with $S_i^1$) such that in each  
$S_i^j$ there exists a {\it unique} genuine fundamental matrix 
\begin{eqnarray}\label{fundsect}
Y_i^j(z)=F_i^j(z)(z-a_i)^{\Lambda_i}e^{Q_i(z)}
\end{eqnarray}
of the system (\ref{syst}) whose factor $F_i^j(z)$ has the asymptotic expansion $\widehat F_i(z)$ in $S_i^j$
(see \cite[Th. 21.13, Prop. 21.17]{IY}). In every intersection $S_i^j\cap S_i^{j+1}$ the fundamental matrices
$Y_i^j$ and $Y_i^{j+1}$ necessarily differ by a constant invertible matrix: 
$$
Y_i^{j+1}(z)=Y_i^j(z)C_i^j, \qquad C_i^j\in{\rm GL}(p,{\mathbb C}),
$$ 
and it is understood that the logarithmic term $(z-a_i)^{\Lambda_i}$ is analytically continued from $S_i^1$ to 
$S_i^2$, from $S_i^2$ to $S_i^3$, ..., from $S_i^{N_i}$ to $S_i^1$, so that 
\begin{eqnarray}\label{circle}
Y_i^1(z)e^{2\pi{\bf i}\Lambda_i}=Y_i^{N_i}(z)C_i^{N_i} \quad\mbox{in}\quad S_i^{N_i}\cap S_i^1.
\end{eqnarray}
The matrices $C_i^1,\ldots,C_i^{N_i}$ are called {\it Stokes matrices} of the system (\ref{syst}) at the 
non-resonant irregular singular point $a_i$. They satisfy the relation
\begin{eqnarray}\label{formmonodr}
e^{2\pi{\bf i}\Lambda_i}=M_i^1\,C_i^1\ldots C_i^{N_i},
\end{eqnarray}
where $M_i^1$ is the monodromy matrix of $Y_i^1$ at the point $a_i$. Indeed, the fundamental matrix $Y_i^1$ can
be continued from $S_i^1$ into $S_i^2$ as $Y_i^2(C_i^1)^{-1}$, since $Y_i^1=Y_i^2(C_i^1)^{-1}$ in 
$S_i^1\cap S_i^2$. Further it is continued from $S_i^2$ into $S_i^3$ as $Y_i^3(C_i^1C_i^2)^{-1}$, etc. Finally, 
in $S_i^{N_i}$ it becomes equal to $Y_i^{N_i}(C_i^1\ldots C_i^{N_i-1})^{-1}$. Then it comes back into $S_i^1$ as
$Y_i^1\,e^{2\pi{\bf i}\Lambda_i}(C_i^1\ldots C_i^{N_i})^{-1}$ according to (\ref{circle}), whence the relation 
(\ref{formmonodr}) follows. It is also known that all the eigenvalues of any Stokes matrix are equal to $1$, that is,
the Stokes matrices are unipotent (see \cite[Prop. 21.19]{IY} or \cite[Th. 15.2]{Wa}).

\section{Linear differential systems and meromorphic connections on holomorphic vector bundles}

Let us recall some notions concerning holomorphic vector bundles and meromorphic connections in
a context of linear differential equations. Here we mainly follow \cite[Ch. 3]{IY} or \cite{GP} 
(see also \cite{BMM}). 

In an analytic interpretation, a holomorphic bundle $E$ of rank $p$ 
over the Riemann sphere is defined by a {\it cocycle} $\{g_{\alpha\beta}(z)\}$, 
that is, a collection of holomorphic matrix functions corresponding to a
covering $\{U_{\alpha}\}$ of the Riemann sphere:
$$
g_{\alpha\beta}: U_{\alpha}\cap U_{\beta}\longrightarrow{\rm GL}
(p,\mathbb C), \qquad U_{\alpha}\cap U_{\beta}\ne\varnothing.
$$
These functions satisfy the conditions
$$
g_{\alpha\beta}=g_{\beta\alpha}^{-1}, \qquad
g_{\alpha\beta}g_{\beta\gamma}g_{\gamma\alpha}=I \quad (\mbox{for
} U_{\alpha}\cap U_{\beta}\cap U_{\gamma}\ne\varnothing).
$$
Two holomorphically equivalent cocycles $\{g_{\alpha\beta}(z)\}$,
$\{g'_{\alpha\beta}(z)\}$ define the same bundle. Equivalence
of cocycles means that there exists a set $\{h_{\alpha}(z)\}$ of
holomorphic matrix functions $h_{\alpha}: U_{\alpha}\longrightarrow
{\rm GL}(p, {\mathbb C})$ such that
\begin{eqnarray}\label{equiv}
h_{\alpha}(z)g_{\alpha\beta}(z)=g'_{\alpha\beta}(z)h_{\beta}(z). 
\end{eqnarray}
A {\it section} $s$ of the bundle $E$ is determined by a set
$\{s_{\alpha}(z)\}$ of vector functions $s_{\alpha}:
U_{\alpha}\longrightarrow{\mathbb C}^p$ that satisfy the
conditions $s_{\alpha}(z)=g_{\alpha\beta}(z) s_{\beta}(z)$ in
intersections $U_{\alpha}\cap U_{\beta}\ne\varnothing$ .

A {\it meromorphic connection} $\nabla$ on the holomorphic vector
bundle $E$ is determined by a set $\{\omega_{\alpha}\}$ of matrix
meromorphic differential 1-forms that are defined in the corresponding
neighbourhoods $U_{\alpha}$ and satisfy gluing conditions
\begin{eqnarray}\label{glue}
\omega_{\alpha}=(dg_{\alpha\beta})g_{\alpha\beta}^{-1}+
g_{\alpha\beta}\omega_{\beta}g_{\alpha\beta}^{-1}\qquad (\mbox{for
} U_{\alpha}\cap U_{\beta}\ne\varnothing).
\end{eqnarray}
Under a transition to an equivalent cocycle $\{g'_{\alpha\beta}\}$
connected with the initial one by the relations (\ref{equiv}), the
1-forms $\omega_{\alpha}$ of the connection $\nabla$ are
transformed into the corresponding 1-forms
\begin{eqnarray}\label{formequiv}
\omega'_{\alpha}=(dh_{\alpha})h_{\alpha}^{-1}+
h_{\alpha}\omega_{\alpha}h_{\alpha}^{-1}.
\end{eqnarray}
Conversely, the existence of holomorphic matrix functions
$h_{\alpha}: U_{\alpha}\longrightarrow{\rm GL}(p, {\mathbb C})$
such that the matrix 1-forms $\omega_{\alpha}$ and
$\omega'_{\alpha}$ (satisfying the conditions (\ref{glue}) for
$g_{\alpha\beta}$ and $g'_{\alpha\beta}$ respectively) are
connected by the relation (\ref{formequiv}) in $U_{\alpha}$, 
indicates the equivalence of the cocycles $\{g_{\alpha\beta}\}$ and
$\{g'_{\alpha\beta}\}$ (one may assume that the intersections 
$U_{\alpha}\cap U_{\beta}$ do not contain singular points of the
connection). 

Vector functions $s_{\alpha}(z)$ satisfying linear differential 
equations  $ds_{\alpha}=\omega_{\alpha} s_{\alpha}$ in the 
corresponding $U_{\alpha}$, by virtue of the conditions(\ref{glue}) 
can be chosen so that a set $\{s_{\alpha}(z)\}$ determines a section 
of the bundle $E$, which is called {\it horizontal} with respect to 
the connection $\nabla$. Thus horizontal sections of a holomorphic 
vector bundle with a meromorphic connection are determined by 
solutions of local linear differential systems. The {\it monodromy of 
a connection} (the monodromy group) characterizes ramification of 
horizontal sections under their analytic continuation along loops in
$\overline{\mathbb C}$ not containing singular points of the
connection 1-forms and is defined similarly to the monodromy group
of the system (\ref{syst}). A connection may be called {\it Fuchsian} 
({\it logarithmic}), {\it regular} or {\it irregular} depending on
the type of the singular points of its 1-forms (as singular points of 
linear differential systems). 

If a bundle is holomorphically trivial (all matrices of the
cocycle can be taken as the identity matrices), then by virtue of the
conditions (\ref{glue}) the matrix 1-forms of a connection coincide 
on non-empty intersections $U_{\alpha}\cap U_{\beta}$. Hence horizontal 
sections of such a bundle are solutions of a global linear differential
system defined on the whole Riemann sphere. Conversely, the linear system 
(\ref{syst}) determines a meromorphic connection on the holomorphically 
trivial vector bundle of rank $p$ over $\overline{\mathbb C}$. It is 
understood that such a bundle has the standard definition by the cocycle 
that consists of the identity matrices while the connection is defined by 
the matrix 1-form $B(z)dz$ of coefficients of the system. But for us it will 
be more convenient to use the following equivalent coordinate description 
(a construction already appearing in \cite{BMM}).

At first we consider a covering $\{U_{\alpha}\}$ of the punctured Riemann 
sphere $\overline{\mathbb C}\setminus\{a_1,\ldots,a_n\}$ by simply connected 
neighbourhoods. Then on the corresponding non-empty intersections 
$U_{\alpha}\cap U_{\beta}$ one defines the matrix functions of a cocycle, 
$g'_{\alpha\beta}(z)\equiv{\rm const}$, which are expressed in terms 
of the monodromy matrices $M_1,\ldots,M_n$ of the system (\ref{syst}) 
via the operations of multiplication and taking the inverse (see \cite{GP}). 
In this case the matrix differential 1-forms $\omega'_{\alpha}$ 
defining a connection are equal to zero. Further the covering
$\{U_{\alpha}\}$ is complemented by small neighbourhoods $O_i$ of
the singular points $a_i$ of the system, thus we obtain the
covering of the Riemann sphere $\overline{\mathbb C}$. To
non-empty intersections $O_i\cap U_{\alpha}$ there correspond
matrix functions $g'_{i\alpha}(z)=Y_i(z)$ of the cocycle, where
$Y_i(z)$ is a germ of a fundamental matrix of the system whose
monodromy matrix at the point $a_i$ is equal to $M_i$ (so, for
analytic continuations of the chosen germ to non-empty
intersections $O_i\cap U_{\alpha}\cap U_{\beta}$ the cocycle
relations $g_{i\alpha}g_{\alpha\beta}=g_{i\beta}$ hold). The matrix
differential 1-forms $\omega'_i$ determining the connection in
the neighbourhoods $O_i$ coincide with the 1-form $B(z)dz$ of
coefficients of the system. To prove holomorphic equivalence of 
the cocycle $\{g'_{\alpha\beta}, g'_{i\alpha}\}$ to the identity 
cocycle it is sufficient to check existence of holomorphic matrix 
functions
$$
h_{\alpha}: U_{\alpha}\longrightarrow{\rm GL}(p, {\mathbb C}),
\qquad h_i: O_i\longrightarrow{\rm GL}(p, {\mathbb C}),
$$
such that
\begin{eqnarray}\label{formequiv2}
\omega'_{\alpha}=(dh_{\alpha})h_{\alpha}^{-1}+
h_{\alpha}\omega_{\alpha}h_{\alpha}^{-1}, \qquad
\omega'_i=(dh_i)h_i^{-1}+ h_i\omega_ih_i^{-1}.
\end{eqnarray}
Since we have $\omega_{\alpha}=B(z)dz$ and $\omega'_{\alpha}=0$
for all $\alpha$, the first equation in (\ref{formequiv2}) is
rewritten as a linear system
$$
d(h^{-1}_{\alpha})=(B(z)dz)h^{-1}_{\alpha},
$$
which has a holomorphic solution $h^{-1}_{\alpha}: U_{\alpha}
\longrightarrow{\rm GL}(p,{\mathbb C})$ since the 1-form
$B(z)dz$ is holomorphic in a simply connected neighbourhood $U_{\alpha}$. 
The second equation in (\ref{formequiv2}) has a holomorphic solution 
$h_i(z)\equiv I$, as $\omega_i=\omega'_i=B(z)dz$.
\medskip

One says that a bundle $E$ has a subbundle $E'\subset E$ of rank
$k<p$ that is {\it stabilized} by a connection $\nabla$, if the
pair $(E, \nabla)$ admits a coordinate description
$\{g_{\alpha\beta}\}$, $\{\omega_{\alpha}\}$ of the following
blocked upper-triangular form:
$$
g_{\alpha\beta}= \left(\begin{array}{cc}
       g_{\alpha\beta}^1 & * \\
       0 & g_{\alpha\beta}^2 \\
      \end{array}\right), \qquad
\omega_{\alpha}= \left(\begin{array}{cc}
       \omega_{\alpha}^1 & * \\
       0 & \omega_{\alpha}^2 \\
      \end{array}\right),
$$
where $g_{\alpha\beta}^1$ and $\omega_{\alpha}^1$ are blocks of
size $k\times k$ (then the cocycle $\{g_{\alpha\beta}^1\}$ defines
the subbundle $E'$ and the 1-forms $\omega_{\alpha}^1$ define the
restriction $\nabla'$ of the connection $\nabla$ to the subbundle
$E'$).
\medskip

{\bf Example 1.} Consider a system (\ref{syst}) whose monodromy matrices 
$M_1,\ldots,M_n$ are of the same blocked upper-triangular form, and the 
corresponding holomorphically trivial vector bundle $E$ with the meromorphic
connection $\nabla$. Suppose that in a neighbourhood of each singular point $a_i$ 
of the system there exist a fundamental matrix $Y_i(z)$ whose monodromy matrix 
is $M_i$, and a holomorphically invertible matrix $\Gamma_i(z)$ such that 
$\Gamma_i\,Y_i$ is a blocked upper-triangular matrix (with respect to the
blocked upper-triangular form of the matrix $M_i$). Let us show that
to the common invariant subspace of the monodromy matrices there corresponds
a vector subbundle $E'\subset E$ that is stabilized by the connection $\nabla$.

We use the above coordinate description of the bundle and connection
with the cocycle $\{g'_{\alpha\beta}, g'_{i\alpha}\}$ and set 
$\{\omega'_{\alpha}, \omega'_i\}$ of matrix 1-forms. The matrices 
$g'_{\alpha\beta}$ are already blocked upper-triangular since the monodromy
matrices $M_1,\ldots,M_n$ are (and $\omega'_{\alpha}=0$), while the matrices  
$g'_{i\alpha}=Y_i$ can be transformed to such a form, $\Gamma_i\,g'_{i\alpha}=
\Gamma_i\,Y_i$. Thus changing the matrices $g'_{i\alpha}$ onto 
$\Gamma_i\,g'_{i\alpha}$ and matrix 1-forms $\omega'_i$ onto
$$
\Gamma_i\,\omega'_i\,\Gamma^{-1}_i+(d\Gamma_i)\Gamma^{-1}_i,
$$
we pass to the holomorphically equivalent coordinate description with the cocycle
matrices and connection matrix 1-forms having the same blocked upper-triangular
form. 
\medskip

The following auxiliary lemma points to a certain block structure
of a linear differential system in the case when the corresponding 
holomorphically trivial vector bundle with the meromorphic connection 
has a holomorphically trivial subbundle that is stabilized by the 
connection.
\medskip

{\bf Lemma 1.} {\it If a holomorphically trivial vector bundle $E$
of rank $p$ over $\overline{\mathbb C}$ endowed with a meromorphic
connection $\nabla$ has a holomorphically trivial subbundle
$E'\subset E$ of rank $k$ that is stabilized by the connection, then
the corresponding linear system $(\ref{syst})$ is reduced to a blocked
upper-triangular form via a constant gauge transformation 
$\tilde y(z)=Cy(z)$, $C\in{\rm GL}(p,\mathbb C)$. That is, 
$$
CB(z)C^{-1}=\left(\begin{array}{cc} B'(z) & * \\
                            0 & *
                  \end{array}\right),
$$
where $B'(z)$ is a block of size $k\times k$.}
\medskip

{\bf Proof.} Let $\{s_1,\ldots,s_p\}$ be a basis of global holomorphic
sections of the bundle $E$ (which are linear independent at each
point $z\in\overline{\mathbb C}$) such that the 1-form of the
connection $\nabla$ in this basis is the 1-form $B(z)dz$ of
coefficients of the linear system. Consider also a basis
$\{s'_1,\ldots,s'_p\}$ of global holomorphic sections of the
bundle $E$ such that $s'_1,\ldots,s'_k$ are sections of the
subbundle $E'$, $(s'_1,\ldots,s'_p)=(s_1,\ldots,s_p)C^{-1}$,
$C\in{\rm GL}(p,{\mathbb C})$.

Now choose a basis $\{h_1,\ldots,h_p\}$ of sections of the bundle
$E$ such that these are horizontal with respect to the connection
$\nabla$ and $h_1,\ldots,h_k$ are sections of the subbundle $E'$
(this is possible since $E'$ is stabilized by the connection). Let
$Y(z)$ be a fundamental matrix of the system whose columns are 
the coordinates of the sections $h_1,\ldots,h_p$ in the basis 
$\{s_1,\ldots,s_p\}$. Then
$$
\widetilde Y(z)=CY(z)=\left(\begin{array}{cc} k\times k & * \\
                                            0 & *
                      \end{array}\right)
$$
is a blocked upper-triangular matrix, since its columns are the
coordinates of the sections $h_1,\ldots,h_p$ in the basis 
$\{s'_1,\ldots, s'_p\}$. Consequently, the transformation 
$\tilde y(z)=Cy(z)$ reduces the initial system to a blocked 
upper-triangular form. {\hfill $\Box$}
\medskip

The degree $\deg E$ (which is an integer) of a holomorphic
vector bundle $E$ endowed with a meromorphic connection $\nabla$ 
may be defined as the sum
$$
\deg E=\sum_{i=1}^n{\rm res}_{a_i}\,{\rm tr}\,\omega_i
$$
of the residues of local differential 1-forms ${\rm tr}\,\omega_i$ over all
singular points of the connection (the notation ''tr'' is for the trace), 
where $\omega_i$ is the local matrix differential 1-form of the connection 
$\nabla$ in a neighbourhood of its singular point $a_i$. Later when calculating 
the degree of a bundle we will apply the Liouville formula ${\rm tr}\,\omega_i=
d\ln\det Y_i$, where $Y_i$ is a fundamental matrix of the local
linear differential system $dy=\omega_i\,y$.

\section{Solvability of regular singular systems with small exponents}

Consider a system (\ref{syst}) with {\it regular} singular points
$a_1,\ldots,a_n$ of Poincar\'e rank $r_1,\ldots,r_n$ respectively.
If the real part of the exponents of this system is sufficiently small,
then the following necessary condition, of solvability by
quadratures, holds.
\medskip

{\bf Theorem 1.} {\it Let for some $k\in\{1,\ldots,p-1\}$ the exponents
$\beta_i^j$ of the regular singular system  $(\ref{syst})$ satisfy the condition
\begin{eqnarray}\label{ineq1}
{\rm Re}\,\beta_i^j>-1/nk, \qquad i=1,\ldots,n, \quad
j=1,\ldots,p,
\end{eqnarray}
and, for each pair $\beta_i^j\not\equiv\beta_i^l$ 
$({\rm mod}\,\mathbb Z)$, $i=1,\ldots,n$, one of the conditions
\begin{eqnarray}\label{ineq2}
{\rm Re}\,\beta_i^j-{\rm Re}\,\beta_i^l\not\in{\mathbb Q} 
\qquad{\rm or}\qquad {\rm Im}\,\beta_i^j\ne{\rm Im}\,\beta_i^l.
\end{eqnarray}
Then the solvability of the system $(\ref{syst})$ by quadratures implies 
the existence of a constant matrix $C\in{\rm GL}(p, \mathbb C)$ such that the 
matrix $CB(z)C^{-1}$ is of the following blocked form:
$$
CB(z)C^{-1}=\left(\begin{array}{cc} B'(z) & * \\
                                        0 & *
                  \end{array}\right),
$$
where $B'(z)$ is an upper-triangular matrix of size $k\times k$.}
\medskip

{\bf Remark 1.} Though the inequalities (\ref{ineq1}) restrict the real parts
of the exponents from below, together with the estimates (\ref{fuchsineq}) 
they provide boundedness from above.

{\bf Remark 2.} The sum of the Poincar\'e ranks of a regular singular system whose 
exponents satisfy the condition (\ref{ineq1}) is indeed restricted 
because of the Fuchs inequalities (\ref{fuchsineq}), namely, $\sum_{i=1}^nr_i<p/k$.
\medskip

The proof of the theorem is based on two auxiliary lemmas.
\medskip

{\bf Lemma 2.} {\it Let the exponents $\beta_i^j$ of the regular singular system 
$(\ref{syst})$ satisfy the condition $(\ref{ineq1})$. If the monodromy matrices of
this system are upper-triangular, then there is a constant matrix 
$C\in{\rm GL}(p,\mathbb C)$ such that the matrix $CB(z)C^{-1}$ has the form as 
in Theorem 1.}
\medskip

{\bf Proof.} We use a geometric interpretation (exposed in the
previous section) according to which to the regular singular system
(\ref{syst}) there corresponds the holomorphically trivial vector
bundle $E$ of rank $p$ over the Riemann sphere endowed with the
meromorphic connection $\nabla$ with the regular singular points $a_1,\ldots,a_n$. 

Since the monodromy matrices $M_1,\ldots,M_n$ of the system are 
upper-triangular there exists, as shown in Example 1, a flag 
$E^1\subset E^2\subset\ldots\subset E^p=E$ of subbundles of rank 
$1,2,\ldots,p$ respectively that are stabilized by the connection 
$\nabla$. Indeed, a fundamental matrix $Y$ determining the monodromy
matrices $M_1,\ldots,M_n$ of the system can be represented near each
singular point $a_i$ in the form 
$$
Y(z)=T_i(z)(z-a_i)^{E_i}, \qquad E_i=(1/2\pi{\bf i})\ln M_i,
$$
where $T_i$ is a meromorphic matrix at the point $a_i$. This matrix can be factored
as $T_i=V_i\,P_i$, with a holomorphically invertible matrix $V_i$ at
$a_i$ and an upper-triangular matrix $P_i$ which is a polynomial in
$(z-a_i)^{\pm1}$ (see, for example, \cite[Lemma 1]{Go}). Thus the matrix
 $V^{-1}_i\,Y$ is upper-triangular.

Let us estimate the degree of each subbundle $E^j$, $j\leqslant k$, noting that
in a neighbourhood of each singular point $a_i$ the initial system is transformed
via a holomorphically invertible gauge transformation in a system with a fundamental
matrix of the form
$$
Y_i(z)=\left(\begin{array}{cc} U'_i(z) & * \\
                                     0 & *
             \end{array}\right)(z-a_i)^
       {\Bigl(\begin{array}{cc} A'_i & 0 \\
                                   0 & *
              \end{array}\Bigr)}(z-a_i)^
       {\Bigl(\begin{array}{cc} E'_i & * \\
                                   0 & *
        \end{array}\Bigr)}
$$
such that the matrix $Y'_i(z)=U'_i(z)(z-a_i)^{A'_i}(z-a_i)^{E'_i}$
is a Levelt fundamental matrix for a linear system of $j$ equations with the
regular singular point $a_i$. The matrix 1-form of coefficients of this system
in a neighbourhood of $a_i$ is a local 1-form of the restriction
$\nabla^j$ of the connection $\nabla$ to the subbundle $E^j$, and the exponents
$\tilde\beta_i^1,\ldots,\tilde\beta_i^j$ of this system (the eigenvalues of the
matrix $A'_i+E'_i$) form a subset of the exponents of the initial 
system at $a_i$. Therefore, 
$$
{\rm Re}\,\tilde\beta_i^l>-1/nk, \qquad l=1,\ldots,j.
$$ 
The degree of the holomorphically trivial vector bundle $E^p$ 
is equal to zero, and for $j\leqslant k$ one has:
\begin{eqnarray*}
\deg E^j & = & \sum_{i=1}^n{\rm res}_{a_i}d\ln\det Y'_i(z)=
\sum_{i=1}^n{\rm ord}_{a_i}\det U'_i(z)+\sum_{i=1}^n{\rm tr}
(\Lambda'_i+E'_i)= \\
 & = & \sum_{i=1}^n{\rm ord}_{a_i}\det U'_i(z)+\sum_{i=1}^n
\sum_{l=1}^j{\rm Re}\,\tilde\beta_i^l>-j/k\geqslant-1.
\end{eqnarray*}
As the degree of a subbundle of a holomorphically trivial vector bundle is 
non-positive, one has $\deg E^j=0$, hence all the subbundles 
$E^1\subset\ldots\subset E^k$ are holomorphically trivial (a subbundle of
a holomorphically trivial vector bundle is holomorphically trivial, if its
degree is equal to zero, see \cite[Lemma 19.16]{IY}). Now the assertion of the lemma
follows from Lemma 1. {\hfill $\Box$}
\medskip

A matrix $A$ will be called $N$-{\it resonant}, if there are two eigenvalues 
$\lambda_1\ne\lambda_2$ such that $\lambda_1^N=\lambda_2^N$, that is,
$$
|\lambda_1|=|\lambda_2|,\qquad \arg\lambda_1-\arg\lambda_2=
\frac{2\pi}N\,j, \quad j\in\{1,2,\ldots,N-1\}.
$$

Let a group $M\subset{\rm GL}(p,{\mathbb C})$ be generated by matrices
$M_1,\ldots,M_n$. If these matrices are sufficiently close to the identity  
(in the Euclide topology) then the existence of a solvable normal subgroup of finite
index in $M$ implies their triangularity, see Theorem 2.7 \cite[Ch. 6]{Kh}. 
According to the remark following this theorem, the requirement of proximity of
the matrices $M_i$ to the identity can be weakened as follows.
\medskip

{\bf Lemma 3.} {\it There is a number $N=N(p)$ such that if the matrices  
$M_1,\ldots,M_n$ are not $N$-resonant, then the existence of a solvable normal 
subgroup of finite index in $M$ implies their triangularity.}
\medskip

{\bf Proof of Theorem 1.} From the theorem assumptions it follows that the identity
component $G^0$ of the differential Galois group $G$ of the system (\ref{syst}) is solvable, 
hence $G^0$ is a solvable normal subgroup of $G$ of finite index. Then the
monodromy group $M$ of this system also has a solvable normal subgroup of 
finite index, namely $M\cap G^0$. 

As follows from the definition of the exponents $\beta_i^j$ of a linear differential
system at a regular singular point $a_i$, these are connected with the eigenvalues 
$\mu_i^j$ of the monodromy matrix $M_i$ by the relation
$$
\mu_i^j=\exp(2\pi\,{\bf i}\,\beta_i^j).
$$
Therefore,
$$
\mu_i^j=\exp(2\pi\,{\bf i}({\rm Re}\,\beta_i^j+{\bf i}\,{\rm Im}\,\beta_i^j))= 
e^{-2\pi\,{\rm Im}\,\beta_i^j}(\cos(2\pi\,{\rm Re}\,\beta_i^j)+
{\bf i}\,\sin(2\pi\,{\rm Re}\,\beta_i^j)),
$$
and for any $N$ the matrices $M_i$ are non $N$-resonant by the conditions 
(\ref{ineq2}) on the numbers $\beta_i^j$. Now the assertion of the theorem follows from 
Lemma 2 and Lemma 3. {\hfill $\Box$}
\medskip

As a consequence of Theorem 1 we obtain the following refinement of the 
Ilyashenko--Khovansky theorem on solvability by quadratures of a Fuchsian system with small 
residue matrices.
\medskip

{\bf Corollary 1.} {\it Let the eigenvalues $\beta_i^j$ of the residue matrices $B_i$ of the
Fuchsian system $(\ref{fuchs})$ satisfy the condition
\begin{eqnarray}\label{ineq3}
{\rm Re}\,\beta_i^j>-\frac1{n(p-1)}, \qquad i=1,\ldots,n, \quad
j=1,\ldots,p,
\end{eqnarray}
and, for each pair $\beta_i^j\not\equiv\beta_i^l$ 
$({\rm mod}\,\mathbb Z)$, $i=1,\ldots,n$, one of the conditions $(\ref{ineq2})$.
Then the solvability of the Fuchsian system $(\ref{fuchs})$ by quadratures is equivalent to 
the simultaneous triangularity of the matrices $B_i$.}
\medskip

{\bf Proof.} The necessity of simultaneous triangularity is a direct consequence of Theorem 1, 
since the exponents of the Fuchsian system (\ref{fuchs}) at $a_i$ 
coincide with the eigenvalues of the residue matrix $B_i$. Sufficiency follows from
a general fact that any linear differential system with an (upper-) triangular 
coefficient matrix is solvable by quadratures (one should begin with the last equation).
{\hfill $\Box$}
\medskip

{\bf Remark 3.} The inequalities (\ref{ineq3}) restricting the real parts of the
exponents of the Fuchsian system from below also provide their boundedness from above
because of the Fuchs relation $\sum_{i=1}^n\sum_{j=1}^p\beta_i^j=0$ (see 
(\ref{fuchsineq})). Namely, 
$$
-\frac1{n(p-1)}<{\rm Re}\,\beta_i^j<\frac{np-1}{n(p-1)}.
$$
In particular, the integer parts $\varphi_i^j$ of the numbers ${\rm Re}\,\beta_i^j$ 
for such a system have to belong to the set $\{-1, 0, 1\}$. 
\medskip

{\bf Remark 4.} If each residue matrix $B_i$ of the Fuchsian system (\ref{fuchs})
only has one eigenvalue $\beta_i$, then the solvability of this system by quadratures
is also equivalent to the simultaneous triangularity of the matrices $B_i$ (not 
depending on values of ${\rm Re}\,\beta_i$). Indeed, in this case each monodromy
matrix $M_i$ only has one eigenvalue $\mu_i=e^{2\pi{\bf i}\beta_i}$, hence is not
$N$-resonant. Then the solvability implies the simultaneous triangularity of the monodromy 
matrices and existence of a flag $E^1\subset E^2\subset\ldots\subset E^p=E$ of subbundles
of the holomorphically trivial vector bundle $E$ that are stabilized by the
logarithmic connection $\nabla$ (corresponding to the Fuchsian system). Since
$\deg E=\sum_{i=1}^np\beta_i=0$, the degree $\sum_{i=1}^nj\beta_i$ of each subbundle
$E^j$ is zero and all these subbundles are holomorphically trivial.
\medskip 

It is natural to expect that for a general Fuchsian system (with no restrictions
on the exponents) solvability by quadratures not necessarily implies simultaneous 
triangularity of the residue matrices. This is indeed illustrated by the following
example of A.\,Bolibrukh.
\medskip

{\bf Example 2} (A.\,Bolibrukh \cite[Prop. 5.1.1]{Bo}). There exist points 
$a_1, a_2, a_3, a_4$ on the Riemann sphere and a Fuchsian system with singularities 
at these points, whose monodromy matrices are
\begin{eqnarray*}
M_1=\left(\begin{array}{rrrrrrr} 1 &  1 &  1 & -1 & 0 & 0 & -1 \\
                                 0 & -1 & -1 &  0 & 0 & 0 &  1 \\
												         0 &  0 &  1 &  1 & 2 & 2 &  2 \\
																 0 &  0 &  0 &  1 & 1 & 0 &  1 \\
																 0 &  0 &  0 &  0 & 1 & 1 &  1 \\
																 0 &  0 &  0 &  0 & 0 & 1 &  0 \\
																 0 &  0 &  0 &  0 & 0 & 0 & -1
					\end{array}\right), \quad
M_2=\left(\begin{array}{rrrrrrr} 1 &  1 &  0 &  1 &  1 &  1 &  0 \\
                                 0 & -1 &  1 &  1 & -1 &  1 & -1 \\
												         0 &  0 & -1 & -1 &  1 & -1 &  0 \\
																 0 &  0 &  0 &  1 &  1 &  1 &  0 \\
																 0 &  0 &  0 &  0 & -1 & -1 &  1 \\
																 0 &  0 &  0 &  0 &  0 &  1 &  0 \\
																 0 &  0 &  0 &  0 &  0 &  0 & -1
					\end{array}\right),\\
M_3=\left(\begin{array}{rrrrrrr} 1 &  0 &  1 &  0 & -1 &  0 &  0 \\
                                 0 & -1 & -1 &  1 & -1 &  1 &  2 \\
												         0 &  0 &  1 &  1 & -1 &  1 &  2 \\
																 0 &  0 &  0 & -1 &  1 & -1 & -2 \\
																 0 &  0 &  0 &  0 & -1 &  1 &  1 \\
																 0 &  0 &  0 &  0 &  0 &  1 &  1 \\
																 0 &  0 &  0 &  0 &  0 &  0 & -1
					\end{array}\right), \quad
M_4=\left(\begin{array}{rrrrrrr} 1 &  0 &  1 & -1 & 1 & 1 &  0 \\
                                 0 & -1 &  1 & -2 & 0 & 0 &  0 \\
												         0 &  0 & -1 &  1 & 0 & 0 &  0 \\
																 0 &  0 &  0 & -1 & 1 & 0 &  1 \\
																 0 &  0 &  0 &  0 & 1 & 1 &  0 \\
																 0 &  0 &  0 &  0 & 0 & 1 &  1 \\
																 0 &  0 &  0 &  0 & 0 & 0 & -1
					\end{array}\right),
\end{eqnarray*}
whereas the coefficient matrix of this system is {\it not} upper-triangular. 
Moreover, this system cannot be transformed in an upper-triangular 
form neither via a constant gauge transformation nor even a 
meromorphic (rational) one preserving the singular points 
$a_1, a_2, a_3, a_4$ and the orders of their poles. Thus the residue 
matrices of this Fuchsian system are not triangular in any basis
though the system is solvable by quadratures, since its monodromy
group generated by the triangular matrices is solvable.

We notice that Corollary 1 does not apply to this example as its exponents 
cannot satisfy the conditions (\ref{ineq3}). Indeed, for any exponent 
$\beta_i^j=\varphi_i^j+\rho_i^j$ one has
$$
\rho_i^j=\frac1{2\pi\bf i}\ln\mu_i^j, \qquad \mu_i^j\in\{-1,1\},
$$ 
hence ${\rm Re}\,\rho_i^j$ is equal to $0$ or $1/2$. The inequalities  
(\ref{ineq3}) imply (for $n=4$, $p=7$)
$$
{\rm Re}\,\beta_i^j>-\frac1{24},
$$
hence $\varphi_i^j$ is equal to $0$ or $1$ (see Remark 3). Therefore,
the sum of the exponents over all singular points is positive, which
contradicts the Fuchs relation.

\section{A criterion of solvability for a non-resonant irregular system with small
formal exponents}

Consider the system (\ref{syst}) with {\it non-resonant} irregular
singular points $a_1,\ldots,a_n$ of Poincar\'e rank $r_1,\ldots,r_n$ 
respectively. If the real part of the formal exponents of this system
is sufficiently small, then the following criterion of solvability by
quadratures holds.
\medskip

{\bf Theorem 2.} {\it Let at each singular point $a_i$ the formal exponents 
$\lambda_i^j$ of the irregular system $(\ref{syst})$ be} pairwise distinct  
{\it and satisfy the condition 
$$
{\rm Re}\,\lambda_i^j>-\frac1{n(p-1)}, 
$$
and, for every pair $(\lambda_i^j, \lambda_i^l)$, one of the conditions 
$(\ref{ineq2})$. Then this system is solvable by quadratures if and
only if there is a constant matrix $C\in{\rm GL}(p, \mathbb C)$ such
that $CB(z)C^{-1}$ is upper-triangular.}
\medskip

{\bf Proof.} As in Corollary 1, sufficiency does not require a special 
proof. Let us prove necessity.

Consider a fundamental matrix $Y$ of the system (\ref{syst}) and the representation
of the differential Galois group $G$ with respect to this matrix. The connection matrices 
between $Y$ and the fundamental matrices $Y_1^1,\ldots,Y_n^1$ (the latter were defined at 
the end of Section 2) denote by $P_1,\ldots,P_n$ respectively: 
$$
Y(z)=Y_i^1(z)P_i, \qquad i=1,\ldots,n.
$$ 
Then, as we noted earlier, the monodromy matrices $M_i=P_i^{-1}M_i^1P_i$ ($i=1,\ldots,n$) 
with respect to $Y$ belong to the group $G$. Moreover, as follows from Ramis' theorem 
\cite{Ra}\footnote{It is difficult to find the original proof of it, but there exist various 
variants of this theorem and comments to it proposed by the other authors (see \cite[Ths. 1, 6]{IKh2}, 
\cite{Ko}, \cite[Th. 2.3.11]{Mi1}+\cite[Prop. 1.3]{Mi2}).}, 
the corresponding Stokes matrices 
$$
\widetilde C_i^1=P_i^{-1}C_i^1P_i,\quad\ldots,\quad
\widetilde C_i^{N_i}=P_i^{-1}C_i^{N_i}P_i\quad (i=1,\ldots,n) 
$$
also belong to $G$. Therefore, by the relation
$$
e^{2\pi{\bf i}\widetilde\Lambda_i}=M_i\widetilde C_i^1\ldots\widetilde C_i^{N_i}, 
\qquad\widetilde\Lambda_i=P_i^{-1}\Lambda_iP_i,
$$ 
(see (\ref{formmonodr})), the group $G$ also contains the matrices 
$e^{2\pi{\bf i}\widetilde\Lambda_i}$ of {\it formal monodromy}.

Denote by $\widehat M$ the group $\bigl<e^{2\pi{\bf i}\widetilde\Lambda_i},\widetilde C_i^1,\ldots,
\widetilde C_i^{N_i}\bigr>_{i=1}^n$ generated by the matrices of formal monodromy and Stokes
matrices over all singular points. As follows from the condition of the theorem, the group $G$ possesses the solvable
normal subgroup $G^0$ of finite index. Hence the subgroup $\widehat M\subset G$ possesses the
solvable normal subgroup of finite index, $\widehat M\cap G^0$. Since for any $N$ the matrices 
generating the group $\widehat M$ are non $N$-resonant, according to Lemma 3 they are simultaniously reduced 
to an upper-triangular form by conjugating to some non-degenerated matrix $\widetilde C$ (non-resonance of the 
formal monodromy matrices $e^{2\pi{\bf i}\widetilde\Lambda_i}$ follows from the conditions on the formal
exponents, the eigenvalues of the matrices $\widetilde\Lambda_i$, and is proved as in Theorem 1; for
the Stokes matrices it follows from their unipotence). We may assume that they are already upper-triangular (otherwise
we would consider the fundamental matrix $Y\widetilde C$ instead of $Y$). As follows from \cite[Ch. VIII, \S1]{Ga}, 
the relation
$$
\Lambda_i=P_i\,\widetilde\Lambda_i\,P_i^{-1},
$$
where $\widetilde\Lambda_i$ is an upper-triangular matrix and $\Lambda_i$ is a diagonal matrix whose diagonal 
elements are pairwise distinct, implies that the matrix $P_i$ writes $P_i=D_i\,R_i$, where $R_i$ is an upper-triangular 
matrix (the conjugation $R_i\,\widetilde\Lambda_i\,R_i^{-1}$ annulates all the off-diagonal elements of the matrix 
$\widetilde\Lambda_i$) and $D_i$ is a permutation matrix for $\Lambda_i$ (that is, the conjugation $D_i^{-1}\,\Lambda_i\,D_i$ 
permutes the diagonal elements of the matrix $\Lambda_i$).

In a neighbourhood of each $a_i$ we pass from the set of fundamental matrices $Y_i^1,\ldots,Y_i^{N_i}$, 
which correspond to the sectors $S_i^1,\ldots,S_i^{N_i}$, to the fundamental matrices  
$$
\widetilde Y_i^j(z)=Y_i^j(z)P_i \qquad (\mbox{in particular, }\widetilde Y_i^1=Y)
$$ 
connected to each other in the intersections $S_i^j\cap S_i^{j+1}$ by the relations 
$$
\widetilde Y_i^{j+1}(z)=\widetilde Y_i^j(z)\widetilde C_i^j.
$$
From the decomposition of $P_i$ above and from (\ref{fundsect}) it follows that the matrices 
$\widetilde Y_i^j$ can be written
$$
\widetilde Y_i^j(z)=F_i^j(z)(z-a_i)^{\Lambda_i}e^{Q_i(z)}D_i\,R_i=F_i^j(z)D_i(z-a_i)^{\Lambda'_i}e^{Q'_i(z)}R_i,
$$
where 
$$
\Lambda'_i=D_i^{-1}\,\Lambda_i\,D_i, \qquad Q'_i(z)=D_i^{-1}\,Q_i(z)\,D_i
$$ 
are diagonal matrices obtained from the corresponding matrices $\Lambda_i$, $Q_i(z)$ by a permutation of the diagonal 
elements. Therefore, in the intersections $S_i^j\cap S_i^{j+1}$ we have the relations
$$
F_i^{j+1}(z)D_i(z-a_i)^{\Lambda'_i}e^{Q'_i(z)}R_i=F_i^j(z)D_i(z-a_i)^{\Lambda'_i}e^{Q'_i(z)}R_i\,\widetilde C_i^j.
$$ 

Thus in the sectors $S_i^1,\ldots,S_i^{N_i}$, which form a covering of a punctured neighbourhood of $a_i$, 
there are holomorphically invertible matrices $F_i^1(z)D_i,\ldots,F_i^{N_i}(z)\,D_i$ respectively such that in the
intersections $S_i^j\cap S_i^{j+1}$ their quotients 
$$
\left(F_i^j(z)D_i\right)^{-1}F_i^{j+1}(z)D_i=(z-a_i)^{\Lambda'_i}e^{Q'_i(z)}R_i\,\widetilde C_i^j\,R_i^{-1}e^{-Q'_i(z)}(z-a_i)^{-\Lambda'_i}
$$
are upper-triangular matrices. As each matrix $F_i^j(z)D_i$ has the same asymptotic expansion $\widehat F_i(z)D_i$ in 
the corresponding $S_i^j$, there exists a matrix $\Gamma_i(z)$ holomorphically invertible at the point $a_i$ such that 
all the matrices 
$$
\widetilde F_i^j(z)=\Gamma_i(z)F_i^j(z)D_i, \qquad j=1,\ldots,N_i,
$$ 
are upper-triangular (according to \cite[Prop. 3]{BJL}). In particular,
$$
\Gamma_i(z)Y(z)=\Gamma_i(z)\widetilde Y_i^1(z)=\Gamma_i(z)F_i^1(z)D_i(z-a_i)^{\Lambda'_i}e^{Q'_i(z)}R_i=
\widetilde F_i^1(z)(z-a_i)^{\Lambda'_i}e^{Q'_i(z)}R_i
$$
is an upper-triangular matrix. Hence (see Example 1) one has a flag 
$E^1\subset E^2\subset\ldots\subset E^p=E$ of subbundles of rank 
$1,2,\ldots,p$ respectively that are stabilized by the connection 
$\nabla$. The bounds on the formal exponents imply that these subbundles 
are holomorphically trivial, whence the assertion of the theorem follows. The proof proceeds as for Lemma 2, there are only two 
small differences: the first one is that now we deal with {\it formal} fundamental matrices of subsystems and their {\it formal} 
exponents, which form a subset of the formal exponents of the initial system at each irregular singular point $a_i$; the second 
one is the appearence of an exponential factor in a formal fundamental matrix of a subsystem, but the logarithmic differential 
of such a factor has a zero residue at $a_i$. {\hfill $\Box$}
\medskip


\end{document}